\documentstyle[amssymb,12pt]{amsart}
\newbox\smilebox
\newbox\anchorbox
\newbox\noanchorbox
\newbox\tempbox

\setbox\smilebox=\hbox{$\smile$}

\def\anchor{\hbox{\vtop{
           \hbox to \wd\smilebox{\hfil\vrule width.4pt height7pt depth1pt\hfil}
           \vskip  -11.5truept
           \hbox to \wd\smilebox{\hfil$\smile$\hfil}}}}
\setbox\anchorbox=\anchor
\def\noanchor{\hbox{\vtop{
           \hbox to \wd\anchorbox{\hfil\anchor\hfil}
           \vskip -14truept
           \hbox to \wd\anchorbox{\hfil/\hfil}}}}
\setbox\noanchorbox=\noanchor

\def\fg#1#2#3{\setbox\tempbox=\hbox{$\scriptstyle{#2}$}
\ifnum\wd\anchorbox>\wd\tempbox\dimen255=\wd\anchorbox
\else\dimen255=\wd\tempbox\fi
{#1\,\vtop{\hbox to \dimen255{\hfil\anchor\hfil}
           \vskip -6truept
           \hbox to \dimen255{\hfil$\scriptstyle{#2}$\hfil}}
           \,#3}}

\def\nfg#1#2#3{\setbox\tempbox=\hbox{$\scriptstyle{#2}$}
\ifnum\wd\noanchorbox>\wd\tempbox\dimen255=\wd\noanchorbox
\else\dimen255=\wd\tempbox\fi
{#1\,\vtop{\hbox to \dimen255{\hfil\noanchor\hfil}
           \vskip -6truept
           \hbox to \dimen255{\hfil$\scriptstyle{#2}$\hfil}}
           \,#3}}

\setbox1=\hbox{$\bot$}

\def\north#1#2{#1\,
\hbox{$\bot$\llap {\hbox to\wd1 {\hfil $/$\hfil}}}
\,#2}

\def\nao#1#2#3{#1\  \hbox{\vtop{ 
\baselineskip=4pt
\hbox{$\bot$\llap {\hbox to\wd1 {\hfil $/$\hfil}}
\hskip .05em \llap{\hbox{$^{\scriptscriptstyle{a}}$}}}\hbox{$\scriptstyle
{#2}$}}}\, #3}

\def\bp{\par{\bf Proof.}$\ \ $}

\def\includeE#1{{\lhook\kern-3.5pt\joinrel\smash{
    \mathop{\longrightarrow}\limits^{#1}}}}

\def\efor/{Example~\ref{E4}}

\def\ep{\par\bigskip}

\def\BL/{Baldwin--Lachlan}
\def\Bu/{Buechler}
\def\Hr/{Hrushovski}
\def\lm/{locally modular}
\def\wm/{weakly minimal}
\def\nm/{non--modular}
\def\tt/{totally transcendental}
\def\ss/{superstable}
\def\ud/{unidimensional}
\def\sm/{strongly minimal}

\def\abar{\overline{a}}

\def\bbar{\overline{b}}
\def\cbar{\overline{c}}

\def\gbar{\overline{g}}
\def\hbar{\overline{h}}

\def\xbar{\overline{x}}
\def\ybar{\overline{y}}
\def\zbar{\overline{z}}

\def\tao{\tau}

\def\dom{{\rm dom}}

\def\tp{{\rm tp}}

\def\tr/{trivial}
\def\nt/{non--trivial}
\def\st/{strong type}

\def\TV/{Tarski--Vaught}

\def\sc/{sound construction}
\def\ac/{atomic construction}

\def\fal/{functional}
\def\upl/{unique parallel lines}
\def\chp/{categorical in a higher power}

\def\text#1{\ \hbox{#1}\ }

\def\sbq{\subset}
\def\contains{\supseteq}
\def\<{\langle}
\def\>{\rangle}
\def\conc{\widehat{~~}}
\def\forces{\mathrel{\raise.4ex\hbox{$\scriptstyle \vert$\hskip-.5ex}\vdash}}

\def\K/{${\cal K}$}
\def\KM{{\cal K}}

\def\r{\restriction}
\def\rl{\r L}
\def\KMOD/{$\KM={\rm Mod}(T_1)\rl$}
\def\LK/{$L(\KM)$}
\def\S{{\rm Spec}(\KM)}
\def\al{\alpha}
\def\b{\beta}

\def\e{\emptyset}
\def\g{\gamma}
\def\k{\kappa}
\def\l{\lambda}
\def\o{\omega}

\def\SN/{\hbox{$\S\ne\e$}}

\def\a1{\aleph_1}

\def\PCD/{${\rm PC}_{\a0}$}
\def\PC/{${\rm PC}_{\Delta}$}
\def\ra{\rightarrow}

\def\bg{\beth_\g}

\def\ss{\smallsetminus}

\def\SG/{\hbox{$\S\cap\bg$}}
\def\SNG/{\hbox{$\S\cap\bg\not =\e$}}
\def\cf{{\rm cf}}

\def\T{{\rm KT}}

\def\phi{\varphi}

\def\endproof{\enspace\vrule height6pt width4pt depth0pt\ep} 
\def\sbq{\subseteq}

\def\KK{{\bf K}}

\def\S{{\cal S}}

\def\T{{\cal T}}
\newtheorem{theorem}{Theorem}[section]

\newtheorem{lemma}[theorem]{Lemma}
\newtheorem{prop}[theorem]{Proposition}

\newtheorem{defn}[theorem]{Definition}
\author{M. C. Laskowski}
\thanks{\hbox{Partially supported by  NSF 
Research~Grants~DMS~9704364 and DMS~0071746}}
\address{Department of Mathematics 
\\ University of Maryland \\ College Park, MD 20742}
\email{mcl@@math.umd.edu}
\author{S. Shelah}
\address{Department of Mathematics\\
Hebrew University of Jerusalem
\and Department of Mathematics\\
Rutgers University}
\thanks{The authors thank the U.S.-Israel
Binational Science Foundation for its support of this project.
This is item 687 in Shelah's bibliography.}
\date{\today}
\title[Karp complexity and the independence property]{Karp complexity and  classes
with the independence property}


\begin{document}

\begin{abstract}
A class $\KK$ of structures is {\em controlled\/} if for all cardinals $\l$,
the relation of $L_{\infty,\l}$-equivalence partitions $\KK$ into a set of equivalence
classes (as opposed to a proper class).  We prove that 
no pseudo-elementary class with the independence
property is controlled.
By contrast, there is a pseudo-elementary class with the strict order
property that is controlled (see~\cite{Karp1}).
\end{abstract}

\maketitle

\section{Introduction}
It is well known that the class of models of an unstable theory
is a rather complicated beast.  Perhaps the most familiar statement of
this complexity is that every such theory $T$ has $2^\kappa$ nonisomorphic models for every
$\kappa>|T|$ (see e.g., \cite{Sh:c}).  In fact, much more is true.  
For instance, in \cite{Sh:e} the second author proves that if $\KK$ is an unsuperstable
pseudo-elementary class (for definiteness $\KK$ is the class of $L$-reducts
of an $L'$-theory $T'$) then for every cardinal $\kappa>|T'|$, $\KK$ contains a 
family of $2^\kappa$ pairwise nonembeddable structures, each of size $\kappa$.

Despite these results, our aim is to give some sort of `classification' to
certain
unstable classes, or to prove that no such classification is possible.  Clearly,
because of the results mentioned above, what is meant by a classification in
this context is necessarily very weak.  Following \cite{Karp1},  a class $\KK$
of structures is
{\em controlled\/} if for every cardinal $\kappa$, the relation
of $L_{\infty,\kappa}$-equivalence partitions $\KK$ into a {\em set\/}
of equivalence classes (as opposed to a proper class of classes).
In \cite{Karp1} we show that this notion has a number of equivalences.
In particular, in \cite{Karp1} we prove  the following proposition (see \cite{Barwise}
or \cite{Karp1} for definitions of the undefined notions):

\begin{prop}  \label{oldchar}
The following are notions are equivalent for any
class $\KK$ of structures.

\begin{enumerate}
\item  $\KK$ is controlled;
\item  For any cardinal $\kappa$, there is an ordinal bound on the $L_{\infty,\kappa}$-Scott
heights of the structures in $\KK$;
\item For any cardinal $\kappa$, there is an ordinal bound on the
$\kappa$-Karp complexity of the structures in $\KK$;
\item  For any cardinal $\mu$, there is a cardinal $\kappa$ such that for
any $M\in\KK$, there are at most $\kappa$ distinct $L_{\infty,\mu^+}$-types
of subsets of $M$ of size at most $\mu$ realized in $M$.
\end{enumerate}
\end{prop}

The whole of this paper is devoted to the proof of the following theorem (see
Definition~\ref{basic}). 

\begin{theorem}  \label{big}
No pseudo-elementary class with the independence property is controlled.
\end{theorem}

To place this result in context, recall that in \cite{Sh:c}, the second author proves
that every unstable theory either has the independence property or has the
strict order property.  Paradigms for these theories are the theory of
the random graph and the theory of dense linear order, respectively.
In \cite{Karp1} we prove that the pseudo-elementary class of doubly transitive linear
orders (which is a subclass of the class of dense linear orders) is controlled.
By contrast, it follows immediately from Theorem~\ref{big} that every 
pseudo-elementary subclass of the class of random graphs is uncontrolled.
That is, with respect to the relation of $L_{\infty,\kappa}$-equivalence,
classes of 
reducts of extensions of the theory of the random graph are sizably more complicated
than certain classes of reducts of extensions of the theory of dense linear order.

The history of this paper is rather lengthy.  The statement of Theorem~\ref{big}
was conjectured by the second author almost ten years ago.  From the outset
it was clear that Theorem~\ref{big} should be proved by 
embedding extremely complicated ordered graphs into structures in $\KK$ using
the generalization of the Ehrenfeucht-Mostowski construction given in Theorem~\ref{fund}.
It was also clear (at least to the second author) that the complexity of the
ordered graph should come from a complicated coloring of pairs from a relatively small cardinal
(see Theorem~\ref{color}).  However, the road from these ideas to a formal proof
was not smooth.  There were a great many false attempts by both authors along the way.
The obstruction was not the infinitary combinatorics.  Rather, it was the very finitary
combinatorics that arose from passing from a well-behaved skeleton to its definable
closure that proved difficult.  

In Section 2 we develop three notions that arise in the proof of Theorem~\ref{big}.
The proof of the theorem is contained in Section~3, with many definitions and
easy lemmas relegated to the appendix.  As the results in the appendix are wholly
self-contained, there is no circularity.

\def\des{{\rm des}}
\def\<{\langle}
\def\>{\rangle}
\def\emptyseq{\<\>}
\def\AA{{\bf A}}
\def\BB{{\bf B}}
\def\etabar{\bar \eta}
\def\alphabar{\bar \alpha}
\def\betabar{\bar \beta}
\def\zetahat{\hat \zeta}
\def\HH{{\cal H}}
\def\etabold{{\bar \eta}}
\def\nubold{{\bar \nu}}
\def\restriction{|}
\def\AA{{\bf A}}
\def\BB{{\bf B}}
\def\alphabar{\bar \alpha}
\def\betabar{\bar \beta}
\def\zetahat{\hat \zeta}

\section{The independence property and complicated colorings}

We begin this preliminary section by proving a
fundamental theorem (Theorem~\ref{fund})
about Skolemized theories with the independence property
and discussing its consequences for pseudo-elementary classes.
Following this, 
we discuss many complicated colorings of certain uncountable cardinals.
We close the section with a short discussion of well-founded trees.

\begin{defn}  
{\em
A formula $\phi(\xbar,\ybar)$ has the {\em independence property\/} with respect to
a theory $T$ if for each $n\in\omega$ there is a model $M$ of $T$
and  sequences
$\<\bbar_i:i<n\>$,
$\<\abar_w:w\sbq n\>$ from $M$ 
such that $M\models \phi(\abar_w,\bbar_i)$ if and only if $i\in w$.

A formula $\psi(\zbar_1,\zbar_2)$ {\em codes graphs\/} if for every (symmetric)
graph $(G,R)$ there is a model $M_G$ of $T$ and $\{\cbar_g:g\in G\}$ from $M_G$
such that for all $g,h\in G$, $M_G\models\psi(\cbar_g,\cbar_h)$ if and only if $R(g,h)$.

A theory $T$ has the {\em independence property\/} if some formula $\phi(\xbar,\ybar)$
has the independence property with respect to $T$.
}
\end{defn}

The next lemma tells us that 
if a theory $T$ has the independence property, then there is a formula 
that both codes graphs and has the independence property with respect to $T$.

\begin{lemma}  \label{coded} Let $T$ be any theory.
\begin{enumerate}
\item If $\psi(\zbar_1,\zbar_2)$ codes graphs, then $\psi(\zbar_1,\zbar_2)$ has the independence 
property with respect to $T$.
\item If $\phi(\xbar,\ybar)$ has the independence property with respect to
$T$, then the formula 
$$\psi(\xbar_1\ybar_1,\xbar_2\ybar_2):=\phi(\xbar_1,\ybar_2)\vee \phi(\xbar_2,\ybar_1)$$
codes graphs.
\end{enumerate}
\end{lemma}

\bp  
(1)  Fix $n$ and let $G=\{g_i:i<n\}\cup\{h_w:w\sbq n\}$ be any symmetric graph with $n+2^n$ vertices
that satisfies $R(g_i,h_w)$ holds if and only if $i\in w$.
Let  $\<\bbar_i:i\in n\>$, $\<\abar_w:w\sbq n\>$ be sequences from some model $M_G$ of $T$
that codes $G$.
Then $M_G\models\psi(\abar_w,\bbar_i)$ if and only if $i\in w$.

(2) It suffices to show that every finite graph can be coded, so fix a finite (symmetric)
graph $(G,R)$ where $G=\{g_i:i<n\}$.  For each $i<n$, let $w_i=\{j<n:R(g_i,g_j)\}$.
Choose a model $M$ of $T$ and sequences 
$\<\bbar_i:i<n\>$,  
$\<\abar_w:w\sbq n\>$  
from $M$ exemplifying the independence property for $\phi(\xbar,\ybar)$.
Let $\cbar_i=\abar_{w_i}\bbar_i$ for each $i<n$.  It is easily verified that
$M\models\psi(\cbar_i,\cbar_j)$ if and only if $R(g_i,g_j)$ holds.  \endproof

Although coding graphs is a desirable property in its own right,
its utility for constructing models 
is greatly increased when it is combined with
an appropriate notion of indiscernibility.
With this objective 
in mind, we generalize the construction of Ehrenfeucht and Mostowski (see e.g., \cite{EM})
to admit skeletons that are indexed by structures that are more complicated than
linear orderings.  
We define an
{\em ordered graph\/} to be a structure ${\cal G}=(G,\le,R)$, where
$\le$ is interpreted as a linear order on $G$ and $R$ is a symmetric, irreflexive 
binary relation.

What makes the class of ordered graphs desirable as index structures is the presence of
the Ne\u{s}et\u{r}il-R\"odl theorem.  The version stated below is sufficient for our
purposes, but is less general than the statement in either 
\cite{AbLeo} or \cite{Nesetril}.

\begin{theorem}  \label{NR}
[Ne\u{s}et\u{r}il-R\"odl Theorem]  For every $e,M\in\omega$ and
every finite ordered graph $P$, there is an ordered graph $Q$ such that for any
coloring $F:[Q]^e\ra M$ there is an ordered subgraph $Y\sbq Q$ that is isomorphic to $P$
such that $F(A)=F(B)$ for any $A,B\in [Y]^e$ that are isomorphic as ordered graphs.
\end{theorem}

The proof of the theorem below is virtually identical with the proof of the
classical Ehrenfeucht-Mostowski theorem, with the 
Ne\u{s}et\u{r}il-R\"odl theorem taking the place of Ramsey's theorem.
Recall that a theory $T$ is {\em Skolemized\/} if every substructure of every model
of $T$ is an elementary substructure.

\begin{theorem} \label{fund}
Let $T$ be any Skolemized theory with the independence property 
and suppose that the formula $\phi(\xbar_1,\xbar_2)$
codes graphs.  
For any ordered graph $G$ there is a model $M_G$ of $T$ and $\{\abar_g:g\in G\}$ from $M_G$
such that 
\begin{enumerate}
\item  The universe of $M_G$ is the definable closure of $\{\abar_g:g\in G\}$;
\item  If $f:H_1\ra H_2$ is any ordered graph isomorphism between finite subgraphs of $G$,
then 
$$M_G\models \psi(\abar_g:g\in H_1)\leftrightarrow \psi(\abar_{f(g)}:g\in H_1)$$
for all formulas $\psi$; and 
\item  For all $g,h\in G$, $M_G\models \phi(\abar_g,\abar_h)$ 
if and only if ${\cal G}\models R(g,h)$.
\end{enumerate}
\end{theorem}

\bp  
If we expand the language to $L(G)$ by adding a sequence of  new constant symbols
$\cbar_g$ for every $g\in G$, then Conditions~(2) and (3) can be expressed by sets of
$L(G)$-sentences.  The consistency of these sentences follows immediately from Lemma~\ref{coded},
Theorem~\ref{NR}
and compactness.\endproof

As notation, we call
$\{\abar_g:g\in G\}$ the {\em skeleton\/} of $M_G$.
Next we extend the notion of independence to pseudo-elementary classes.

\begin{defn}  \label{basic}
{\em
Fix a language $L$.  A class $\KK$ of $L$-structures is a 
{\em pseudo-elementary class\/} if there is a language $L'\contains L$ and
an $L'$-theory $T'$ such that $\KK$ is the class of $L$-reducts of models of $T'$.
Such a class has the {\em independence property\/} if some $L$-formula
$\phi(\xbar,\ybar)$ has the independence property with respect to $T'$.
}
\end{defn}

Note that as we can always assume that $T'$ is Skolemized, the conclusions of Theorem~\ref{fund}
apply to any pseudo-elementary class.  The caveat is that in Clause~(1), every element of $M_G$
will be in the $L'$-definable closure of the skeleton, where $L'$ is the language of the Skolemized 
theory.

Our method of proving Theorem~\ref{big} will be to use the theorem
above to produce a family of elements of $\KK$ that code some very complicated 
ordered graphs.  To make this complexity explicit, we discuss some properties of colorings
that were developed by the second author.  See \cite{Sh:g} for a more complete account
of these notions.
As notation, for $x$ a finite subset of a cardinal $\mu$,
let $x^m$ denote the $m^{{\rm th}}$ element of $x$ in 
increasing order.  

\begin{theorem}  \label{color}
Suppose that $\mu=\kappa^{++}$ for any infinite cardinal $\kappa$.
There is a symmetric two-place function $c:\mu\times\mu\ra \mu$ such that
for every $n\in\o$, every collection of $\mu$ disjoint, $n$-element subsets 
$\{x_\al:\al\in\mu\}$ of $\mu$,
and every function $f:n\times n\ra \mu$,
there are $\al<\b<\mu$ such that 
$$c(x^m_\al,x^{m'}_\b)=f(m,m')$$
for all $m,m'<n$.
\end{theorem}

The existence of such a coloring $c$ is called $Pr_0(\mu,\mu,\mu,\aleph_0)$
in both \cite{Sh:g} and \cite{Sh:572}.  The same notion is called $Pr^+(\mu)$
in \cite{Sh:327}.  
Theorem~\ref{color} follows immediately from the results in \cite{Sh:327}
for all uncountable $\kappa$ (since the set $S_{\kappa^+}=\{\al\in\kappa^{++}:
\cf(\al)=\kappa^+\}$ is nonreflecting and stationary).  The case of $\kappa=\aleph_0$
is somewhat special and is proved in \cite{Sh:572} by a separate argument.

%
%
%
%
%
%

We close this section by recalling the definition of a well-founded tree
and proving an easy coloring lemma.

\begin{defn}  \label{otree} 
{\em An {\em $\o$-tree\/}
$\T$ is a downward closed subset of $^{<\o}\l$ for some ordinal $\l$.
We call $\T$ {\em well-founded\/} if it does not have an infinite branch.
For a well-founded 
tree $\T$ and $\eta\in \T$, the {\em depth of $\T$ above $\eta$\/}, ${\rm dp}_{\T}(\eta)$
is defined inductively by
$${\rm dp}_{\T}(\eta)=\left\{\begin{array}{ll}
\sup\{{\rm dp}_{\T}(\nu)+1:\eta\lessdot\nu\}&\mbox{if $\eta$ has a successor}\\
                           0&\mbox{otherwise.}
\end{array}\right. $$

and the depth of $\T$, ${\rm dp}(\T)={\rm dp}_{\T}(\langle\rangle)$.
}
\end{defn}

The most insightful example is that for any ordinal $\delta$, the tree $(\des(\delta),\lessdot)$ 
consisting of
all descending sequences of ordinals less than $\delta$ ordered by initial segment
has depth $\delta$.

\begin{lemma}  \label{tree} 
If $\T\sbq {^{<\o}\l}$ is well-founded and has depth $\kappa^+$, then
any coloring $f:\T\ra \k$, there is a sequence $\<a_n:n\in\omega\>$
of elements from $\T$ such that  $lg(a_n)=n$ and 
$f(a_m|_n)=f(a_n)$ for all $n\le m<\omega$.
\end{lemma}

\bp
For each $n\in\o$ we will find a subset $X_n\sbq \k^+$ of size $\k^+$
and a function $g_n:X_n\ra \T\cap {^n\l}$ such that $X_{n+1}\sbq X_n$,
every element of $g_{n+1}(X_{n+1})$ is a successor of an element of $g_n(X_n)$,
${\rm dp}_\T(g_n(\al))\ge\al$, and $f|_{g_n(X_n)}$ is constant.

To begin, let $X_0=\k^+$ and let $g_0:X_0\ra \{\langle\rangle\}$.
Given $X_n$ and $g_n$ satisfying our demands, we
define  $X_{n+1}$ and $g_{n+1}:X_{n+1}\ra \T\cap {^{n+1}\l}$ as follows:
For $\al\in X_n$, let $\b$ be the least element of $X_n$ greater than $\al$.
As ${\rm dp}_{\T}(g_n(\b))\ge\b$, we can define $g_{n+1}(\al)$ to be a successor of
$g_n(\b)$ of depth at least $\al$.  Since $X_n$ has size $\k^+$, let $X_{n+1}$
be a subset of $X_n$ of size $\k^+$ such that $f|_{g_{n+1}(X_{n+1})}$ is monochromatic.

Now for each $n\in\o$,
simply take $a_n=g_n(\b_n)$, where $\b_n$ is the least element of $X_n$.
\endproof

\section{Proof of Theorem~\ref{big}}

Fix any pseudo-elementary class $\KK$ with the independence property.  
For definiteness, suppose that $L\sbq L'$ are languages and $T'$ is an $L'$-theory such
that $\KK$ is the class of $L$-reducts of models of $T'$.
Without loss, we may assume that $T'$ is Skolemized.  Let $\mu=|T'|^{++}$.
Fix an $L$-formula $\phi(x_1,x_2)$ that codes graphs (see Lemma~\ref{coded}).
For notational simplicity we assume that $lg(x_1)=lg(x_2)=1$.

Now assume by way of contradiction that $\KK$ is controlled.  
It follows from Proposition~\ref{oldchar}(4) that there is a cardinal $\kappa$ such that
for any $M\in\KK$ there are fewer than $\kappa$ distinct $L_{\infty,\mu^+}$-types of subsets
of size at most $\mu$ in $M$.  {\bf Fix,
for the whole of the paper, such a $\kappa$ and put $\delta:=\kappa^+$.}

Our strategy for proving Theorem~\ref{big} is to define one specific
structure $M_\delta\in\KK$.  This $M_\delta$ is
constructed using Theorem~\ref{fund} and is the reduct
the Skolem Hull of  an ordered graph $I_\delta$.
The ordering on $I_\delta$ is a well-order, 
but the edge relation on $I_\delta$ is
extremely complicated as it codes a coloring $c$ given by Theorem~\ref{color}.
The definition of $I_\delta$ and the construction of $M_\delta$ 
are completed in the paragraph
following Definition~\ref{R0}.
Following our construction of $M_\delta$ we use the bound on the number
of $L_{\infty,\mu^+}$-types to form an $\omega$-sequence 
$\<\BB_n:n\in\omega\>$ of $\mu$-sequences of pairs of elements from $M_\delta$
that are reasonably coherent.  Then, by combining several of the results
from the Appendix with properties of the coloring $c$ 
we establish three claims 
whose statements follow Definition~\ref{int}.
These claims collectively imply the existence of
an infinite, descending sequence of ordinals below $\delta$.
This contradiction demonstrates that the class $\KK$ is not controlled.

\begin{defn}  \label{des} 
{\em  The expression $\des(\delta)$ denotes the set of
all strictly decreasing sequences of elements from $\delta$.
The set of all finite sequences from $\des(\delta)$ is denoted by
$\des{^{<\omega}}(\delta)$.
}
\end{defn}

Every element of $\des(\delta)$ is clearly a finite sequence.
It is an easy exercise to show that $\des(\delta)$ 
is a well-ordering with respect to  the lexicographic order  $<_{lex}$.
As noted in the remarks following Definition~\ref{otree}, the
$\omega$-tree $(\des(\delta),\lessdot)$ has depth $\delta=\kappa^+$.

\begin{defn}
{\em    A function $g:\zeta\ra\des{^{<\omega}}(\delta)$
is {\em uniform\/} if $\omega\le\zeta\le\mu$;
$lg(g(\al))=lg(g(\b))$ for all $\al,\b\in\zeta$; and (letting $lg(g)$ denote this common length)
for all $i<lg(g)$, the sequences $\<g(\b)(i):\b\in\zeta\>$ have
constant length and are either constant or $<_{lex}$-strictly increasing.
Let $lg(g($--$)(i))$ denote the length of $g(\b)(i)$ for some (every) $\b\in \zeta$.
If the sequence 
$\<g(\b)(i):\b\in\zeta\>$ 
is constant we let $g_i$ denote its common value.
Let ${\cal U}$ denote the set of all uniform functions.
}
\end{defn}

\begin{defn}
{\em   The universe of $I_\delta$ is the set of all
$t=\<\zeta^t,\eta^t,g^t,p^t\>$, where
\begin{enumerate}
\item $\zeta^t\in\mu$;
\item $\eta^t\in\des(\delta)$;
\item $g^t:\zeta^t\ra\des{^{<\omega}}(\delta)$ is a uniform function; and
\item $p^t\in\{0,1\}$.
\end{enumerate}
}
\end{defn}

We well-order $I_\delta$ as follows.  First, choose any well-ordering $<_{\cal U}$ on the set 
${\cal U}$ of uniform functions.  Then, define the ordering on $I_\delta$
to be lexicographic i.e., $s<_{I_{\delta}} t$ if and only if either $\zeta^s<\zeta^t$; or
$\zeta^s=\zeta^t$ and $\eta^s<_{\rm lex}\eta^t$; or $\zeta^s=\zeta^t$ and $\eta^s=\eta^t$
and $g^s<_{\cal U} g^t$; or
$\zeta^s=\zeta^t$ and $\eta^s=\eta^t$
and  $g^s= g^t$ and $p^s<p^t$.

In order to define the edge relation on $I_\delta$ we require   some preparatory
definitions.

\begin{defn} 
{\em  Two uniform functions $g$ and $h$ (possibly with different domains)
have the {\em same shape\/}
if the following four conditions hold:
\begin{enumerate}
\item  $lg(g)=lg(h)$;
\item  For each $i<lg(g)$, $lg(g($--$)(i))=lg(h($--$)(i))$;
\item  For each $i<lg(g)$,
the sequence
$\<g(\b)(i):\b\in\dom(g)\>$ is constant 
if and only if
$\<h(\b)(i):\b\in\dom(h)\>)$ is constant;
\item  For all  $i,j<lg(g)$ such that $\<g(\b)(i):\b\in\dom(g)\>$
and $\<g(\b)(j):\b\in\dom(g)\>$ are both constant,
$g_i=g_j\Leftrightarrow h_i=h_j$ and
$g_i\lessdot g_j\Leftrightarrow h_i\lessdot h_j$.
\end{enumerate}
}
\end{defn}

\begin{defn} 
{\em  Two pairs $(s,t), (s',t')\in (I_\delta)^2$
have the {\em same type\/} if the following conditions hold:
\begin{enumerate}
\item $\zeta^s<\zeta^t\Leftrightarrow \zeta^{s'}<\zeta^{t'}$ and
$\zeta^s>\zeta^t\Leftrightarrow \zeta^{s'}>\zeta^{t'}$;
\item $lg(\eta^t)=lg(\eta^{t'})$;
\item $p^s=p^{s'}$ and $p^t=p^{t'}$;
\item The uniform functions $g^s$ and $g^{s'}$ have the same shape;
\item For all $i<lg(g^s)$, $g^s(\zeta^t)(i)=\eta^t\Leftrightarrow g^{s'}(\zeta^{t'})(i)=\eta^{t'}$ and
$g^s(\zeta^t)(i)\lessdot\eta^t\Leftrightarrow g^{s'}(\zeta^{t'})(i)\lessdot\eta^{t'}$.
\end{enumerate}

Evidently, having the same type induces an equivalence relation on pairs from $I_\delta$
with countably many classes.  We let $\tp(s,t)$ denote the class of pairs that have
the same type as $(s,t)$ and let ${\cal E}$ denote the set of equivalence classes.
Let $\HH$ denote any countable collection of (total) functions from ${\cal E}$ to $\{0,1\}$
such that for any partial function $f:{\cal E}\ra\{0,1\}$ whose domain is finite
there is an $h\in \HH$ extending $f$.
}
\end{defn}

Using Theorem~\ref{color} choose a symmetric, binary function
$$c:\mu\times\mu\ra\HH$$
such that for every $k\in\omega$, for every collection of $\mu$ disjoint, $k$-element
subsets $\{x_\al:\al\in\mu\}$ of $\mu$, and for every function $f:k\times k\ra\HH$,
there are $\al<\b<\mu$ such that $c(x_\al^m,x_\b^{m'})=f(m,m')$
for all $m,m'<k$.  (Here, $x^m_\al$ denotes the $m^{\rm th}$ element of $x_\al$.)

We are now able to complete our description of the ordered graph $I_\delta$
by defining the edge relation $R(x,y)$ on $I_\delta$.

\begin{defn}  \label{R0}
{\em 
For $s,t\in I_\delta$, $R_0(s,t)$ holds if and only if the following
conditions are satisfied:
\begin{enumerate}
\item  $\zeta^s>\zeta^t$;
\item  $lg(\eta^s)<lg(\eta^t)$;
\item  $p^s=1$; $p^t=0$; and
\item  $c(\zeta^s,\zeta^t)[\tp(s,t)]=1$.
\end{enumerate}
We say $R(s,t)$ holds if and only if $R_0(s,t)$ or $R_0(t,s)$ holds.
}
\end{defn}

Let $M_\delta'\models T'$ be a model satisfying the conclusions of Theorem~\ref{fund}
with respect to the ordered graph $(I_\delta,\le_{I_\delta},R)$
defined above.  To ease notation, we identify
the $I_\delta$ with the skeleton $\{a_g:g\in I_\delta\}$ of $M_\delta'$.
In particular, every element of $M'_\delta$ is an $L'$-term
applied to a finite sequence from $I_\delta$.
Let $M_\delta\in\KK$ be the
$L$-reduct of $M_\delta'$.

As notation, let $g^\al_{\emptyseq}$ denote the 
function whose domain is $\alpha$
and $g(\beta)=\emptyseq$ for all $\beta\in \alpha$.
For $\nu\in\des(\delta)$, let $A_{\nu,\al}\in M_\delta^2$ denote the pair of
elements $\<{(\al,\nu,g^{\al}_{\emptyseq},0)},
{(\al,\nu,g^{\al}_{\emptyseq},1)}\>$ from $I_\delta$
(recall that we are identifying $I_\delta$ with the skeleton)
and let $\AA_\nu$ denote the sequence $\<A_{\nu,\al}:\al\in\mu\>$.

As the number of $L_{\infty,\mu^+}$-types of subsets of $M_\delta$ of size at most
$\mu$ is bounded 
by $\kappa$ and $lg(\AA_\nu)\le\mu$ for all $\nu\in\des(\delta)$,
there is a function $f:\des(\delta)\rightarrow\kappa$ such
that $f(\nu)=f(\nu')$ if and only if $lg(\nu)=lg(\nu')$ and 
$\tp_{\infty,\mu^+}(\<\AA_{\nu|l}:l\le lg(\nu)\>)=
\tp_{\infty,\mu^+}(\<\AA_{\nu'|l}:l\le lg(\nu)\>)$.
Since the depth of the $\omega$-tree $(\des(\delta),\lessdot)$
is $\delta=\kappa^+$, it follows from Lemma~\ref{tree} applied to
this function $f$ that there is a sequence
$\<\nu_n^*:n\in\omega\>$ of elements from $\des(\delta)$
such that for all $n\in\omega$, $lg(\nu^*_n)=n$ 
and the sequences $\<\AA_{\nu^*_n|l}:l\le n\>$ and
$\<\AA_{\nu^*_m|l}:l\le n\>$  have the same $L_{\infty,\mu^+}$-type in $M_\delta$
for all $m\ge n$.

Thus, one can construct by induction on $n$
an $\omega$-sequence $\<\BB_n:n\in\omega\>$ in $M_\delta$
such that:
\begin{itemize}
\item  Each $\BB_n$ is a sequence $\<B_{n,\alpha}:\al\in\mu\>$, where each $B_{n,\al}$
is a pair of elements from $M_{\delta}$;
\item $\BB_0=\AA_{\<\>}$; and
\item The sequences $\<\BB_l:l\le n\>$ and  $\<\AA_{\nu^*_n|l}:l\le n\>$ 
have the same $L_{\infty,\mu^+}$-type for every $n\in\omega$.
\end{itemize}

Fix sequences $\<\nu^*_n:n\in\omega\>$ and $\<\BB_n:n\in\omega\>$ satisfying the properties 
described above.
As notation, we write $a_\alpha$ for the element
$(\al,\emptyseq,g^\al_{\emptyseq},1)\in I_\delta$ (i.e., the second coordinate of 
$A_{\emptyseq,\al}$).  For $n>0$ we write
$a_{n,\al}$ for $(\al,\nu^*_n,g^\al_{\emptyseq},0)$ (the
first coordinate of $A_{\nu^*_n,\al}$) and write $b_{n,\al}$ for the first
coordinate of $B_{n,\alpha}$.  We let $\Gamma_n=\tp(a_\al,a_{n,\b})$ for all
$\al>\b$ from $\mu$.  So, for example, when $\al>\b$ then
$$M_{\delta}\models \phi(a_\al,b_{n,\b})\Leftrightarrow
M_{\delta}\models \phi(a_\al,a_{n,\b})\Leftrightarrow c(\al,\b)[\Gamma_n]=1.$$

Next we use results from the Appendix to obtain subsequences of the 
sequences $\<a_\al:\al\in\mu\>$ and $\<b_{n,\al}:\al\in\mu\>$ with
desirable regularity properties.  
It may be helpful to the reader to skip ahead to the Appendix
at this point in order to become familiar with the definitions
therein.  Specifically,
iterating Lemma~\ref{cleanpairlemma} yields a
descending sequence $Y_1\contains Y_2\contains\dots$ of stationary
subsets of $\mu$ such that for all $n>0$ the sequences $\<a_\al:\al\in Y_n\>$
and $\<b_{n,\al}:\al\in Y_n\>$ form a clean pair (see Definition~\ref{cleanpair}).
As notation, for each $n>0$
fix a number $m(n)$, an $L'$-term $\tao_n$, and for each $l<m(n)$,
tidy sequences $\<b^l_{n,\al}:\al\in Y_n\>$
of elements from the skeleton $I_\delta$
such that for all $\al\in Y_n$,
$$b_{n,\al}=\tao_n(b^l_{n,\al}:l<m(n)).$$
Let
$\zeta^l_{n,\al}$, $\eta^l_{n,\al}$, $g^l_{n,\al}$, and $p^l_{n,\al}$
denote the four components of $b^l_{n,\al}$.

Fix an $n>0$.  Let $Y^*_n=\{(\al,\b)\in (Y_n)^2:\al>\b\}$.
It follows from the remarks following Definition~\ref{clean}
that $$b^0_{n,\b}<_{I_\delta} b^1_{n,\b}<_{I_\delta}\dots <_{I_\delta} b^{m(n)-1}_{n,\b}<_{I_\delta}
a_\al$$
for all $(\al,\b)\in Y^*_n$ and that
$R(b^l_{n,\b},b^{l'}_{n,\b})\leftrightarrow 
R(b^l_{n,\b'},b^{l'}_{n,\b'})$ for all $\b,\b'\in Y_n$ and for all $l,l'<m(n)$.
Thus, 
the only freedom we have in determining whether $\phi(a_\al,b_{n,\b})$ holds or fails
for various $(\al,\b)\in Y_n^*$ is whether or not $R(a_\al,b^l_{n,\b})$
holds or fails for various $l<m(n)$.
Accordingly, we call a subset $Z\sbq m(n)$ {\em true for $n$\/} if
$$M_{\delta}\models \phi(y,\tao_n(x^l:l<m(n)))$$
for all $<_{I_{\delta}}$-increasing sequences $x^0,\dots,
x^{m(n)-1}, y$ from $I_\delta$
such that $R(x^l,x^{l'})$ holds if and only if 
$R(b^l_{n,\b},b^{l'}_{n,\b})$ holds for $\b\in Y_n$ and $R(y,x^l)$ holds if and only if $l\in Z$.
A subset $Z\sbq m(n)$ is {\em false for $n$\/} if it is not true for $n$.

We call an index $l\in m(n)$ {\em $n$-constant\/} if $\zeta^l_\b=\zeta^l_{\b'}$ for all $\b,\b'\in Y_n$.
Let $\b^*_l$ denote this common value.  
As $\<a_\al:\al\in Y_n\>$ and $\<b_{n,\al}:\al\in Y_n\>$ form a clean pair,
it follows that for every $n$-constant $l$,
the values of both $c(\al,\b^*_l)$  and $\tp(a_\al,b^l_{n,\b})$ are constant 
for all $(\al,\b)\in Y_n^*$.  Thus, for all $n$-constant $l$,
$$R(a_\al,b^l_{n,\b})\leftrightarrow R(a_{\al'},b^l_{n,\b})$$
for all $(\al,\b),(\al',\b')\in Y^*_n$.  Let $P_n$ denote the set of all
$n$-constant $l$'s such that $R(a_{\al},b^l_{n,\b})$ holds for all 
$(\al,\b)\in Y^*_n$.

Switching our attention to the non-constants, let $J_n$ denote the
set of non-constant $l\in m(n)$. Let 
$$V_n=\{l\in J_n:\zeta^l_\b=\b\ \hbox{and}\ \tp(a_\al,b_{n,\b}^l)=\Gamma_n\
\hbox{for all (some)}\ (\al,\b)\in Y_n^*\}.$$
There is a natural equivalence relation $E_n$ on $V_n$ defined by $E_n(l,l')$
if and only if $\eta^l_\b=\eta^{l'}_\b$ for all $\b\in Y_n$.
(It follows from Condition~\ref{clean4} of Definition~\ref{clean} that
whether or not 
$\eta^l_\b=\eta^{l'}_\b$ is independent of $\b$.)
We are now able to state the crucial definition for the argument that follows.

\begin{defn}  \label{int}
{\em 
An $E_n$-class $C$ is {\em $n$-interesting\/}
if there is a union of $E_n$-classes $X\sbq V_n$ such that $P_n\cup X$ is false for
$n$, while $P_n\cup X\cup C$ is true for $n$.
}
\end{defn}

In what follows, we will prove the following three claims.

\medskip\noindent{{\bf Claim 1.}}  For every $n>0$ there is an $n$-interesting 
$E_n$-class $C$.

\medskip\noindent{{\bf Claim 2.}}  For every $n>0$ and for every  $n$-interesting 
$E_n$-class $C$ there is an  $\eta^C\in\des(\delta)$ of length $n$ such that
$\eta^l_\b=\eta^C$ for all $l\in C$ and all $\b\in Y_n$.

\medskip\noindent{{\bf Claim 3.}}  For every $n'>n>0$ and for every  $n$-interesting 
$E_n$-class $C$ there is an $n'$-interesting $E_{n'}$-class $C'$ such that
$\eta^C\lessdot\eta^{C'}$.
\medskip

Clearly, one can deduce a contradiction from the three claims by building an
infinite, descending sequence of ordinals.  Thus, to  complete the 
proof of Theorem~\ref{big} it
suffices to prove the claims.  
The proofs of all three appeal to the complexity
of the coloring $c$.  The first application is direct, but the other two
involve constructing appropriate surrogates to the $a_\al$'s before invoking
the properties of the coloring.

\medskip
{\bf Proof of Claim 1.}  Fix $n>0$.  Let $\alphabar=\{\al\}$, let 
$\betabar=\{\b\}\cup\{
\zeta^l_\b:l\in J_n\}$, and choose $k\ge|\betabar|$.  Note that by
Condition~\ref{clean5} of Definition~\ref{clean}, $\al>\zeta^l_\b$
for all $(\al,\b)\in Y^*_n$.  Let $\pi_1:J_n\ra k$
be the function defined by $\pi_1(l)=t$ if and only if 
$\zeta^l_\b$ is the $t^{\rm th}$ element of $\betabar$.
For each $l\in J_n$ let 
$\Gamma_l=\tp(a_\al,b^l_{n,\b})$ for all $(\al,\b)\in Y_n^*$.  (As $g^{a_\al}$ is the trivial
function and as $\<b_{n,\al}:\al\in Y_n\>$ is clean, it is easily verified that
there is only one such type for each $l\in J_n$.)

Let $h,h':k\times k\ra \HH$ be any functions that satisfy:
\begin{enumerate}
\item $h(0,\pi_1(l))[\Gamma_l]=0$ for all $l\in J_n$;
\item $h(0,0)[\Gamma_n]=0$; and
\item $h'=h$ EXCEPT that $h'(0,0)[\Gamma_n]=1$.
\end{enumerate}

It follows easily from the properties of the coloring $c$ that there is $(\al,\b)\in Y_n^*$
such that $c(\al,\zeta^{\pi_1(l)}_\b)=h(0,\pi_1(l))$ for all $l\in J(n)$.
Fix such a pair $(\al,\b)$ and choose $(\al',\b')\in Y_n^*$ such that
$c(\al',\zeta^{\pi(l)}_{\b'})=h'(0,\pi_1(l))$ for all $l\in J(n)$.
It is readily verified that
$$\{l\in m(n):R(a_\al,b^l_{n,\b})\ \hbox{holds}\}=P_n,$$
while
$$\{l\in m(n):R(a_{\al'},b^l_{n,\b'})\ \hbox{holds}\}=P_n\cup V_n.$$
But, as $c(\al,\b)[\Gamma_n]=0$ and $c(\al',\b')[\Gamma_n]=1$,
$$M_{\delta}\models \neg\phi(a_\al,b_{n,\b})\wedge\phi(a_{\al'},b_{n,\b'}),$$
so $P_n$ is false for $n$, while $P_n\cup V_n$ is true for $n$.

Let $\<C_j:j<s\>$ be an enumeration of the $E_n$-classes of $V_n$.
Choose $j<s$ such that $P_n\cup\bigcup_{i<j} C_i$ is false for $n$,
while $P_n\cup\bigcup_{i\le j} C_i$ is true for $n$.  Then $C_j$ is 
$n$-interesting.\endproof

{\bf Proof of Claim 2.}  Fix $n>0$ and an $n$-interesting
$E_n$-class $C$. Choose $X\sbq V_n$, $X$ a union of $E_n$-classes,
such that $P_n\cup X$ is false for $n$, while $P_n\cup X\cup C$ is true for $n$.

Using Lemma~\ref{total}, choose a stationary subset $W\sbq Y_n$ and a uniform
function $g:\mu\ra\des{^{<\omega}}(\delta)$ that satisfies
$$g(\b)=\<\eta^l_{n,\b}:l<m(n)\> \quad \hbox{for all $\b\in W$}.$$

For all $\al\in W$, let $e_\al$ denote the element $(\al,\emptyseq,g|_\al,1)$
from the skeleton  $I_\delta$ of $M_\delta$.
By applying Lemma~\ref{cleanpairlemma} and possibly shrinking $W$,
we may additionally assume that the sequences $\<e_\al:\al\in W\>$ and 
$\<b_{n,\al}:\al\in W\>$
form a clean pair.
The $e_\al$'s should be thought of as being a surrogate for the $a_\al$'s
that carry just enough data from the $b_{n,\b}$'s.

Let $\Gamma^*_l=\tp(e_\al,b^l_{n,\b})$ for all $\al>\b$ from $W$.
The fact that the values of these types does not depend on our choice of $(\al,\b)$
follows from our choice of the functions $g^{e_{\al}}$ 
and Condition~(\ref{clean4}) of Definition~\ref{clean} applied to $\<b_{n,\al}:\al\in W\>$.
To elaborate, the crucial point is that from our definition of $g^{e_\al}|_W$, relations
such as `$g^{e_\al}(\zeta^l_\b)(l')=\eta^l_\b$' are essentially unary (depending only
on $\b$) when restricted to pairs $\al>\b$ from $W$.
Note that for each $l<m(n)$
the type
$\Gamma_l^*(x,y)$ contains the relation 
\begin{equation}
\eta^y=g^x(\zeta^y)(l).  \label{part}
\end{equation}
Let $\Gamma^*_C=\Gamma^*_l$ for any $l\in C$.  

As $\<a_{n,\al}:\al\in W\>$ realizes the same $L_{\infty,\mu^+}$-type as
$\<b_{n,\al}:\al\in W\>$, we can choose $\<d_\al:\al\in W\>$
such that the sequences
\begin{equation}
\<d_\al:\al\in W\>\conc\<b_{n,\al}:\al\in W\>\
\hbox{and}\
\<e_\al:\al\in W\>\conc\<a_{n,\al}:\al\in W\>. \label{matching}
\end{equation}
have the same $L_{\infty,\mu^+}$-type.

By applying Lemma~\ref{cleanpairlemma} we can find a stationary subset $Z\sbq W$
such that the sequences $\<d_\al:\al\in Z\>$ and $\<a_{n,\al}:\al\in Z\>$ form a
clean pair.  Let $Z^*=\{(\al,\b)\in Z^2:\al>\b\}$.
For each $\al\in Z$ say $$d_\al=\theta(d_\al^r:r<r(d)),$$
where  $\theta$ is an $L'$-term and $\<d^r_\al:r<r(d)\>$ is a strictly $<_{I_\delta}$-increasing
sequence from $I_\delta$.  As notation, let $\zetahat^r_\al$ denote the
$\zeta$-component of $d^r_\al$.  Let $J_d=\{r\in r(d):\zetahat^r_\al$ is not
constant$\}$.  
For each $r<r(d)$, let $\Phi_r=\tp(d^r_\al,a_{n,\b})$ for
any $(\al,\b)\in Z^*$.  Note that since $g^{d_\al}($--$)(r)$ is strictly increasing or constant
for each $\al\in Z$ and $\eta^{a_{n,\b}}=\nu^*_n$ for all $\b\in Z$,
both of the relations 
$$g^{d_\al}(\b)=\eta^{a_{n,\b}}\quad\hbox{and}\quad
g^{d_\al}(\b)\lessdot\eta^{a_{n,\b}}$$
concentrate on tails for all $r<r(d)$ (see Definition~\ref{goodrelation}).  
Thus, it follows from Lemma~\ref{tail}
that by possibly trimming $Z$ further, we may assume that for each $r<r(d)$,
the value of $\Phi_r$ is independent of our choice of 
$(\al,\b)\in Z^*$.

\medskip\noindent{\bf Subclaim.}  There is an $r\in J_d$ such that $\zetahat^r_\al=\al$ 
for $\al\in Z$ and $\Phi_r=\Gamma^*_C$
\smallskip

\bp  Let $\alphabar=\{\al\}\cup\{\zetahat^r_\al:r\in J_d\}$, let
$\betabar=\{\b\}\cup\{\zeta^l_{n,\al}:l\in J_n\}$ and choose $k\ge |\alphabar|,|\betabar|$.
Let $\pi_0:J_d\ra k$ be the function that satisfies $\pi_0(r)=s$ if and only if
$\zetahat^r_\al$ is the $s^{\rm th}$ element of $\alphabar$ and let $\pi_1:J_n\ra k$
be the function that satisfies $\pi_1(l)=t$ if and only if
$\zeta^l_{n,\al}$ is the $t^{\rm th}$ element of $\betabar$.
Since the sequences $\<d_\al:\al\in Z\>$ and $\<b_{n,\al}:\al\in Z\>$ are clean,
the lengths of $\alphabar$, $\betabar$ and the values of $\pi_0$ and $\pi_1$
do not depend on our choice of $\alpha\in Z$.

Now, if the subclaim were false we could find two functions $h,h':k\times k\ra\HH$ that satisfy
the following conditions:
\begin{enumerate}
\item $h(\pi_0(r),\pi_1(l))[\Phi_r]=h'(\pi_0(r),\pi_1(l))[\Phi_r]$ for $r\in J_d$ and 
$l\in J_n$;
\item $h(\pi_0(r),0)[\Gamma^*_l]=h'(\pi_0(r),0)[\Gamma^*_l]=1$ for $r\in J_d$, $l\in X$;
\item $h(\pi_0(r),0)[\Gamma^*_l]=h'(\pi_0(r),0)[\Gamma^*_l]=0$ for $r\in J_d$, $l\in V_n\setminus X
\setminus C$;
\item $h(0,0)[\Gamma^*_C]=0$; $h'(0,0)[\Gamma^*_C]=1$.
\end{enumerate}

\relax From the properties of the coloring $c$, choose $(\al,\b)$ and $(\al',\b')$ from 
$Z^*$ such that
$$c(\zetahat^{\pi_0(r)}_\al,\zeta^{\pi_1(l)}_\b)=h(\pi_0(r),\pi_1(l)) \quad
\!\hbox{and}\quad\!
c(\zetahat^{\pi_0(r)}_{\al'},\zeta^{\pi_1(l)}_{\b'})=h'(\pi_0(r),\pi_1(l))$$
for all $r\in J_d$ and all $l\in J_n$.
Thus,
$$\{l\in m(n):R(e_\al,b^l_{n,\b})\ \hbox{holds}\}=P_n\cup X,$$
which is false for $n$, while
$$\{l\in m(n):R(e_{\al'},b^l_{n,\b'})\ \hbox{holds}\}=P_n\cup X\cup C,$$
which is true for $n$.  Hence
$$M_{\delta}\models \neg\phi(e_\al,b_{n,\b})\wedge\phi(e_{\al'},b_{n,\b'}),$$
so it follows from Equation~(\ref{matching}) that
\begin{equation}
M_{\delta}\models \neg\phi(d_\al,a_{n,\b})\wedge\phi(d_{\al'},a_{n,\b'}).  \label{subclaim}
\end{equation}
However, as $\<d_\al:\al\in Z\>$ and $\<a_{n,\al}:\al\in Z\>$ form a clean pair,
the sequences $\<a_{n,\b}\>\conc\<d^r_\al:r<r(d)\>$ and 
$\<a_{n,\b'}\>\conc\<d^r_{\al'}:r<r(d)\>$ are both $<_{I_\delta}$-strictly increasing.
As well, 
$R(d^r_\al,d^{r'}_\al)\leftrightarrow R(d^r_{\al'},d^{r'}_{\al'})$
holds for all $r,r'<r(d)$ 
by the remark following Definition~\ref{clean}.  Since
$c(\zetahat^r_\al,\b)[\Phi_r]=c(\zetahat^r_{\al'},\b')[\Phi_r]$ for all $r\in J_d$,
$R(d^r_\al,a_{n,\b})\leftrightarrow R(d^r_{\al'},a_{n,\b'})$
holds
for all $r<r(d)$ as well. That is, 
the pairs $(\al,\b)$ and $(\al',\b')$ generate isomorphic
ordered subgraphs of $I_\delta$.  Hence
$$M_{\delta}\models \phi(d_\al,a_{n,\b})\leftrightarrow\phi(d_{\al'},a_{n,\b'}),$$
which contradicts Equation~(\ref{subclaim}).\endproof

To complete the proof of Claim 2 choose any $r<r(d)$
such that $\Phi_r=\Gamma^*_C$ and $\zetahat^r_\al=\al$ for all $\al\in Z$.
As well, fix $\al>\b>\b'$ from $Z$, let $\gbar$ denote the $g$-component from
$d^r_\al$, and choose any $l^*\in C$.
Since $\tp(d^r_\al,a_{n,\b})=\tp(d^r_\al,a_{n,\b'})=\Gamma^*_C$, it follows from
Equation~(\ref{part}) that
$$\gbar(\b)(l^*)=\nu^*_n=\gbar(\b')(l^*).$$
Since $\gbar$ is uniform, the function $\gbar($--$)(l^*)$ must be constant.
As well, this information is part of the shape of $\gbar$.  However, since 
$\tp(e_\al,b^{l^*}_{n,\b})=\tp(d^r_\al,a_{n,\b})$, the function $g^{e_\al}$ has the
same shape as $\gbar$, so the function
$g^{e_\al}($--$)(l^*)$ must be constant as well.  
But the $l^*$-th coordinate of $g^{e_\al}(\b)$ was chosen to be $\eta^{l^*}_{n,\b}$
for all $\beta\in W$.  That is, $\<\eta^{l^*}_{n\b}:\b\in W\>$
is constant.  But, as  the sequence 
$\<\eta^{l^*}_{n,\b}:\b\in Y_n\>$ forms a $\Delta$-system, it too must be constant.
Let $\eta^C$ denote the common value  of $\eta^{l^*}_{n,\b}$. 
That $\eta^l_{n,\b}=\eta^C$
for all $l\in C$ and all $\b\in Y_n$ follows immediately from Condition~(\ref{clean4}) of
Definition~\ref{clean} and the definition of $E_n$.

Finally, since $\tp(a_\al,b^{l}_{n,\b})=\Gamma_n$ for all $l\in V_n$ and all $(\al,\b)\in Y_n$,
$lg(\eta^{l^*}_{n,\b})=n$ as required.\endproof

{\bf Proof of Claim 3.}  Fix $n'>n>0$ and an $n$-interesting
$E_n$-class $C$. By reindexing, we may assume that the index sets
$J_n$ and $J_{n'}$ are disjoint.  Choose $X\sbq V_n$, $X$ a union of $E_n$-classes,
such that $P_n\cup X$ is false for $n$, while $P_n\cup X\cup C$ is true for $n$.

As we are choosing between finitely many possibilities,
by shrinking $Y_{n'}$ further, we may assume that for all $l,l'\in J_n\cup J_{n'}$
the truth values of
the relations 
$$\hbox{`$\eta^l_\al=\eta^{l'}_{\zeta^{l}_\al},$'}\quad
\hbox{`$\eta^l_\al\lessdot\eta^{l'}_{\zeta^{l}_\al}$,'}
\quad\hbox{and}\quad\hbox{`$\eta^{l'}_{\zeta^{l}_\al}\lessdot\eta^l_\al$,'}$$
are invariant among all $\al\in Y_{n'}$.
By analogy with the argument in Claim~2,
use Lemma~\ref{total} to find a stationary subset $W\sbq Y_{n'}$
and a uniform
function $g:\mu\ra\des{^{<\omega}}(\delta)$ that satisfies
$$g(\b)=\<\eta^l_{n,\b}:l<m(n)\>\conc\<\eta^{l'}_{n',\b}:l'<m(n')\>
\quad \hbox{for all $\b\in W$}.$$

For all $\al\in W$, let $e_\al$ denote the element $(\al,\emptyseq,g|_\al,1)$
from the skeleton of $M_\delta$.
(These $e_\al$'s are not the same as in the proof of Claim~2 as the function $g$ is different.)

Let $\Gamma^*_{n,l}=\tp(e_\al,b^l_{n,\b})$ 
and let $\Gamma^*_{n',l'}=\tp(e_\al,b^{l'}_{n',\b})$ 
for all $\al>\b$ from $W$.
As was the case in the proof of Claim~2, the values of $\Gamma^*_{n,l}$ and
$\Gamma_{n',l'}$ do not depend on our choice of $(\al,\b)$.
The verification of this depends on Condition~(\ref{clean4}) of Definition~\ref{clean}
and the further reduction performed above.
Note that for each $l<m(n)$
the type
$\Gamma_{n,l}^*(x,y)$ contains the relation 
`$\eta^y=g^x(\zeta^y)(l)$,' while the type $\Gamma_{n',l'}^*(x,y)$
contains the relation
`$\eta^y=g^x(\zeta^y)(m(n)+l')$,' for all $l'<m(n')$.
As well, note that if $E_n(l_1,l_2)$, then $\Gamma_{n,l_1}^*=\Gamma_{n,l_2}^*$.
Let $\Gamma^*_C=\Gamma^*_l$ for any $l\in C$.  

As $\<a_{n,\al}:\al\in W\>\conc\<a_{n',\al}:\al\in W\>$ realizes the same $L_{\infty,\mu^+}$-type as
$\<b_{n,\al}:\al\in W\>\conc\<b_{n',\al}:\al\in W\>$,
we can choose $\<d_\al:\al\in W\>$ from $M_{\delta}$
such that
$$\<d_\al:\al\in W\>\conc\<b_{n,\al}:\al\in W\>\conc\<b_{n',\al}:\al\in W\>$$
has the same $L_{\infty,\mu^+}$-type as
\begin{equation}
\<e_\al:\al\in W\>\conc\<a_{n,\al}:\al\in W\>\conc\<a_{n',\al}:\al\in W\>
\label{triples}
\end{equation}

Using Lemma~\ref{cleanpairlemma}, choose a stationary subset $Z\sbq W$
such that both pairs of sequences $\{d_\al:\al\in Z\>$, $\{a_{n,\al}:\al\in Z\>$ and
$\{d_\al:\al\in Z\>$, $\{a_{n',\al}:\al\in Z\>$ are clean pairs.
Let $Z^*=\{(\al,\b)\in Z^2:\al>\b\}$.
For each $\al\in Z$ say
$$d_\al=\theta(d_\al^r:r<r(d)),$$
where $\theta$ is an $L'$-term and
$\<d^r_\al:r<r(d)\>$ is a strictly $<_{I_\delta}$-increasing
sequence from $I_\delta$.  As notation, let $\zetahat^r_\al$ denote the
$\zeta$-component of $d^r_\al$.  Let $J_d=\{r\in r(d):\zetahat^r_\al$ is not
constant$\}$.  As in the proof of Claim~2, we can use Lemma~\ref{tail} to shrink
$Z$ so that the values of $\tp(d^r_\al,a_{k,\b})$ is independent of the choice of
$(\al,\b)\in Z^*$ for all $r<r(d)$ and all $k\in\{n,n'\}$.
Let $\Phi_r=\tp(d^r_\al,a_{n,\b})$ for all $(\al,\b)\in Z^*$.

Let $$\alphabar=\{\al\}\cup\{\zetahat^r_\al:r\in J_d\}, 
\betabar=\{\b\}\cup\{\zeta^l_{n,\al}:l\in J_n\}\cup\{\zeta^{l'}_{n',\al}:l'\in J_{n'}\}$$ 
and choose $k\ge |\alphabar|,|\betabar|$.
(Recall that we chose the index sets $J_n$ and $J_{n'}$ to be disjoint.)
Let $\pi_0:J_d\ra k$ be the function that satisfies $\pi_0(r)=s$ if and only if
$\zetahat^r_\al$ is the $s^{\rm th}$ element of $\alphabar$ and let $\pi_1:J_n\cup J_{n'}\ra k$
be the function that satisfies 
$\pi_1(l)=t$ if and only if
$l\in J_n$ and $\zeta^l_{n,\al}$ is the $t^{\rm th}$ element of $\betabar$ {\bf OR}
$l\in J_{n'}$ and 
$\zeta^{l}_{n',\al}$ is the $t^{\rm th}$ element of $\betabar$.
As was the case in the proof of Claim~2, the lengths of $\alphabar$ and $\betabar$ and the 
functions $\pi_0$ and $\pi_1$ do not depend on $\al\in Z$.

%

Suppose that $\Phi(x,y)$ is any type that satisfies $lg(\eta^y)=n$.
We call a type $\Psi$ an {\em extension of $\Phi$\/}
if there are $s,t,t'$ from $I_\delta$ such that $lg(\eta^t)=n$, $lg(\eta^{t'})=n'$,
$\tp(s,t)=\Phi$, $\tp(s,t')=\Psi$, and $\eta^t\lessdot \eta^{t'}$.
Note that any type $\Phi$ has only finitely many extensions.  As well, note that one of the
types $\Gamma^*_{n',l'}$ is an extension of $\Gamma^*_C$, then necessarily $l'\in V_{n'}$.

We call a function $h:k\times k\ra\HH$ {\em closed under $r$-extensions\/}
if $$h(\pi_0(r),0)[\Phi_r]=h(\pi_0(r),0)[\Psi]$$
for all $r\in J_d$ and all of the (finitely many) types $\Psi$ extending $\Phi_r$.

Since $\tp(x,a_{n',\b})$ is an extension of $\tp(x,a_{n,\b})$ for any $\b$ and any
$x$ from $I_\delta$, it follows easily 
that if $h$ is closed under $r$-extensions and some $(\al,\b)\in Z^*$ satisfies
$$c(\zetahat^{\pi_0(r)}_\al,\zeta^{\pi_1(l)}_\b)=h(\pi_0(r),\pi_1(l)) \quad
\!\hbox{and}\quad\!
c(\zetahat^{\pi_0(r)}_{\al'},\zeta^{\pi_1(l)}_{\b'})=h'(\pi_0(r),\pi_1(l))$$
for all $r\in J_d$ and all $l\in J_n\cup J_{n'}$,
then $$R(d^r_\al,a_{n,\b})\leftrightarrow R(d^r_\al,a_{n',\b})\quad\hbox{for all $r<r(d)$},$$
(recall that if $\zetahat^r_\al$ is constant then it follows from cleaning that
$\b>\zetahat^r_\al$ for all $\beta\in Z$, so $R(d^r_\al,a_{n,\b})$ and $R(d^r_\al,a_{n',\b})$
both fail). So, the ordered graph with universe $\{d^r_\al:r<r(d)\}\cup\{a_{n,\b}\}$ is
isomorphic to the ordered graph with universe
$\{d^r_\al:r<r(d)\}\cup\{a_{n',\b}\}$, hence
\begin{equation}
M_{\delta}\models\phi(d_\al,a_{n,\b})\leftrightarrow\phi(d_\al,a_{n',\b}) \label{B}
\end{equation} 
for any such $(\al,\b)\in Z^*$.
Let $$D'=\{l'\in J_{n'}:\pi_1(l')=0\ \hbox{and}\ \Gamma^*_{n',l'}\ \hbox{extends}\ \Gamma^*_C\}$$
and let
$$X'=\{l'\in J_{n'}:\pi_1(l')=0\ \hbox{and}\ \Gamma^*_{n',l'}\ \hbox{extends}\ \Gamma^*_{n,l}\ 
\hbox{for some}\ l\in X\}.$$
Clearly, both $D'$ and $X'$ are subsets of $V_{n'}$ and are unions of $E_{n'}$-classes.

Now fix any function $h:k\times k\ra\HH$ that is closed under $r$-extensions
and satisfies the following conditions:
\begin{enumerate}  
\item For all $l\in J_n$
$$h(0,\pi_1(l))[\Gamma^*_{n,l}]=\left\{ \begin{array}{ll}
                                                  1 & \mbox{if $l\in X$} \\
                                                  0 & \mbox{otherwise;} 
                                                  \end{array} \right.$$
\item For all  $l'\in J_{n'}$
$$h(0,\pi_1(l'))[\Gamma^*_{n',l'}]=\left\{ \begin{array}{ll}
                                                  1 & \mbox{if $l'\in X'$} \\
                                                  0 & \mbox{otherwise;} 
                                                  \end{array} \right.$$
\item For all $r\in J_d$
$$h(\pi_0(r),0)[\Phi_r]=\left\{ \begin{array}{ll}
                                                  1 & \mbox{if $\Phi_r=\Gamma^*_{n,l}$ for some $l\in X$} \\
                                                  0 & \mbox{otherwise.} 
                                                  \end{array} \right.$$
                                           
\end{enumerate}

It is a routine (but somewhat lengthy) exercise to show that there indeed is such a
function $h$.  The key observations are that $X$ and $X'$
are unions of $E_n$ and $E_{n'}$-classes respectively, and that for $k=n$ or $k=n'$,
for all $l_1,l_n\in V_k$,
$$\Gamma^*_{k,l_1}=\Gamma^*_{k,l_2}\quad\hbox{if and only if}\quad E_k(l_1,l_2).$$

Choose any $(\al,\b)\in Z^*$ that satisfies  
$c(\zetahat^{\pi_0(r)}_\al,\zeta^{\pi_1(l)}_\b)=h(\pi_0(r),\pi_1(l))$
for all $r\in J_d$ and all $l\in J_n\cup J_{n'}$.
It follows from Conditions~(1) and (2)
of the constraints on $h$
that
$$\{l\in m(n):R(a_\al,b^l_{n,\b})\ \hbox{holds}\}=P_n\cup X,$$
while
$$\{l'\in m(n'):R(a_\al,b^{l'}_{n',\b})\ \hbox{holds}\}=P_{n'}\cup X'.$$
But $X$ was chosen so that $P_n\cup X$ is false for $n$, hence 
$$M_{\delta}\models \neg\phi(e_\al,b_{n,\b}).$$
It follows from elementarity and the fact that $h$ is closed under $r$-extensions that
$$M_{\delta}\models\neg\phi(d_\al,a_{n,\b})\Rightarrow
M_{\delta}\models\neg\phi(d_\al,a_{n',\b})\Rightarrow
M_{\delta}\models\neg\phi(e_\al,b_{n',\b}),$$
so $P_{n'}\cup X'$ is false for $n'$.

But now, consider the function $h':k\times k\ra\HH$, where $h'=h$
EXCEPT that 
$$h'(\pi_0(r),0)[\Gamma^*_C]=h'(\pi_0(r),0)[\Psi]=1$$
for all types $\Psi$ extending $\Gamma^*_C$.
Note that $h'$ is also closed under $r$-extensions.
Using the properties of the coloring $c$, choose
$(\al',\b')\in Z^*$ such that
$c(\zetahat^{\pi_0(r)}_{\al'},\zeta^{\pi_1(l)}_{\b'})=h'(\pi_0(r),\pi_1(l))$
for all $r\in J_d$ and all $l\in J_n\cup J_{n'}$.
It is easily verified that
$$\{l\in m(n):R(a_{\al'},b^l_{n,\b'})\ \hbox{holds}\}=P_n\cup X\cup C,$$
and
$$\{l'\in m(n'):R(a_{\al'},b^{l'}_{n',\b'})\ \hbox{holds}\}=P_{n'}\cup X'\cup D'.$$
But $M_{\delta}\models \phi(e_{\al'},b_{n,\b'})$.  So, arguing as above,
it follows that
$$M_{\delta}\models \phi(e_{\al'},b_{n,\b'}).$$
Thus, $P_{n'}\cup X'\cup D'$ is true for $n'$.

But now, simply write $D=\{C_0',\dots,C'_{s-1}\}$, where the $C_i$'s are distinct $E_{n'}$-classes.
Thus, there is $j<s$ such that $P_{n'}\cup X'\cup\bigcup_{i<j} C_i$ is false for $n'$,
while $P_n\cup X'\cup\bigcup_{i\le j} C'_i$ is true for $n'$.  In particular,
the class $C'_j$ is ${n'}$-interesting and $\eta^C\lessdot \eta^{C_j'}$ since
$C_j'\sbq D$.\endproof


\appendix

\section{Cleaning Lemmas}  \label{cleaning}
In the appendix we define a number of desirable properties of sequences and show
that if the original sequence was indexed by a stationary subset of $\mu$
(which is regular) then
there is a subsequence that is also indexed by a stationary set that has this desirable 
property.  Many of these properties are unary, which makes  the situation
easy. For example, if every element of the sequence
has one of fewer than $\mu$ colors, then there is a monochromatic stationary subsequence.
It would certainly be desirable to extend this to pairs, i.e.,
if $S\sbq \mu$ is stationary and every pair $(\al,\b)\in S^2$ with $\al>\b$ is given
one of fewer than $\mu$ colors, then one could find a subsequence that is homogeneous
in this sense.  However, for an arbitrary coloring, this would require $\mu$ to be
weakly compact. In fact, the existence of the coloring of pairs
given by Theorem~\ref{color} can be viewed as a strong refutation of 
the existence in general of such a homogeneous set.
However, if we restrict to  relations 
that concentrate on tails (see Definition~\ref{goodrelation}) then
Lemma~\ref{tail} provides us with a
stationary homogeneous subset.

Nothing in this appendix is at all deep.  The arguments simply rely on standard methods
of manipulating clubs and stationary sets, with Fodor's lemma playing a prominent role.
The notation in the appendix is consistent with the body of
the paper.  In particular, the $\mu$, $\delta$, $I_\delta$ and $M_\delta$
that appear in the Appendix are the same entities as in Section~3.

\begin{lemma}  \label{constant} Suppose that $S\sbq \mu$ is stationary and
$f$ is any ordinal-valued function with domain $S$.  Either there is a stationary subset
$S'\sbq S$ such that \hbox{$f\restriction_{S'}$} is constant or there is a stationary
subset $S'\sbq S$ such that $f\restriction_{S'}$ is strictly increasing.
\end{lemma}

\bp  Choose $\delta^*$ least such that there is a stationary $S'\sbq S$
such that $f(\al)<\delta^*$ for all $\al\in S'$.  Without loss, we may assume that
$S'=S$, i.e., $f(\al)<\delta^*$ for all $\al\in S$.
Let $$T=\{\al\in S:f(\al)<f(\b)\ \hbox{for some}\ \b\in S\cap\al\}.$$
We claim that $T$ is not stationary.  Indeed, if $T$ were stationary, then
the function $g:T\ra\mu$ defined by $g(\al)$ is the least $\b\in S$ such that $f(\al)<f(\b)$
would be pressing down.  Thus, by Fodor's lemma there would be 
a stationary $T'\sbq T$ and $\b^*\in S$ such that $g(\al)=\b^*$ for all $\al\in T'$.
But then, $\al\in T'$ would imply $f(\al)<f(\b^*)<\delta^*$, which
contradicts our choice of $\delta^*$.  Thus, $T$ is not stationary.
So by replacing $S$ by $S\setminus T$, we may assume that $f(\al)\ge f(\b)$ for
all $\al<\b$ from $S$.  Let $$U=\{\al\in S:f(\al)=f(\b)\ \hbox{for some}\ \b\in S\cap\al\}.$$
There are now two cases.  If $U$ is stationary then it follows from Fodor's lemma 
that $f$ is constant on some stationary subset of $U$.  On the other hand, 
$f$ is strictly increasing on $S\setminus U$, so if $U$ is non-stationary then the second
clause of the conclusion of the lemma holds.\endproof

\begin{defn}  \label{Deltasystem}
{\em 
For $X\sbq\mu$, a sequence $\etabold=\<\eta_\al:\al\in X\>$
of elements from $\des(\delta)$ forms a {\em $\Delta$-system indexed by $X$\/}
if
\begin{enumerate}
\item  $lg(\eta_\al)=lg(\eta_\b)$ for all $\al,\b\in X$.  This common
value,
called the {\em length of $\etabold$,}
is denoted $lg(\etabold)$;
\item  For each $i<lg(\etabold)$, $\<\eta_\al(i):\al\in X\>$ is either constant
or strictly increasing;
\item For all $i<j<lg(\etabold)$, $\eta_\al(i)\neq\eta_\b(j)$ for all $\al,\b\in X$.
\end{enumerate}
We call $i<lg(\etabold)$ {\em constant\/} if the sequence 
$\<\eta_\al(i):\al\in X\>$ is constant.
}
\end{defn}


\begin{lemma}  \label{Deltaexist}
If $S\sbq\mu$ is stationary, then for any sequence $\<\eta_\al:\al\in S\>$ from $\des(\delta)$
there is a stationary $S'\sbq S$ such that $\<\eta_\al:\al\in S'\>$ is a $\Delta$-system
indexed by $S'$.
\end{lemma}

\bp  The first clause of Definition~\ref{Deltasystem} follows easily from the fact
that the countable union of non-stationary sets is non-stationary and the second clause
follows by iterating Lemma~\ref{constant} finitely often.
To obtain the third clause, assume that the original sequence satisfies the first
two clauses and fix $i<j<lg(\etabold)$.  By the definition of $\des(\delta)$,
$\eta_\al(i)>\eta_\al(j)$ for all $\al\in S$.  If both $i$ and $j$ are constant there is
nothing to do.  If $i$ is constant and $j$ is strictly increasing then necessarily
$\eta_\al(i)>\eta_\b(j)$ for all $\al,\b\in S$ and if $j$ is constant then again
$\eta_\al(i)>\eta_\b(j)$ for all $\al,\b\in S$.  So assume that both sequences
$\<\eta_\al(i):\al\in S\>$ and
$\<\eta_\al(j):\al\in S\>$ are strictly increasing.  It suffices to show
that the set $$T=\{\al\in S:\eta_\al(j)\in\{\eta_\b(i):\b\in S\cap\al\}\}$$
is non-stationary.  However, if $T$ were stationary then for each $\al\in T$,
choose $\b\in S$ least such that $\eta_\b(i)=\eta_\al(j)$.  Since 
$\eta_\b(j)\lessdot\eta_\b(i)=\eta_\al(j)$ and since
$\<\eta_\al(j):\al\in S\>$ is strictly increasing, $\al>\b$.
Thus, Fodor's lemma would give us $\al\neq\al'$ such that $\eta_\al(j)=\eta_{\al'}(j)$,
which contradicts the fact that
$\<\eta_\al(j):\al\in S\>$ is strictly increasing.  \endproof

\begin{defn}  \label{tidy}
{\em 
A sequence $\<s_\al:\al\in X\>$ of elements from $I_\delta$
is {\em tidy\/} if the following conditions hold:
\begin {enumerate}
\item  The sequence $\<\zeta^{s_\al}:\al\in X\>$ is either constant or is
strictly increasing with $\zeta^{s_\al}\ge\al$ for all $\al\in X$;
\item  The sequence $\<\eta^{s_\al}:\al\in X\>$ is a $\Delta$-system indexed by $X$;
\item  The sequence $\<p^{s_\al}:\al\in X\>$ is constant; and
\item  The uniform functions $g^{s_\al}$ and $g^{s_\b}$ have the same shape for all $\al,\b\in X$.
\end{enumerate}
}
\end{defn}

\begin{lemma} \label{tidylemma} If $S\sbq\mu$ is stationary and $\<s_\al:\al\in S\>$ is any sequence
of elements from $I_\delta$, then there is a stationary $S'\sbq S$ such that
the subsequence
$\<s_\al:\al\in S'\>$ is tidy.
\end{lemma}

\bp  The first condition can be obtained by applying Lemma~\ref{constant}
to the sequence $\<\zeta^{s_\al}:\al\in S\>$ to get a subsequence indexed by a stationary
subset $S_1\sbq S$ that is either
constant or strictly increasing.  If the subsequence is strictly increasing,
then  it follows easily from Fodor's lemma that $\{\al\in S_1:\zeta^{s_\al}<\al\}$ is non-stationary
so by trimming $S_1$ further we may assume it is empty.  The second condition follows immediately
from Lemma~\ref{Deltaexist} and the final two conditions can be obtained by noting
that the union of countably many non-stationary subsets of $\mu$ is non-stationary.\endproof


%
%

\begin{defn}  \label{clean}
{\em 
A sequence $\<b_\al:\al\in X\>$ of elements from $M_\delta$ is
{\em clean\/} if there is a term $\tao(x_0,\dots,x_{m-1})$ with $m$ free variables
and sequences $\<s^l_\al:\al\in X\>$ from the skeleton $I_\delta$ for each $l<m$ such
that $$b_\al=\tao(s^0_\al,\dots,s^{m-1}_\al) \ \hbox{for each}\ \al\in X$$
and satisfy the following conditions (as notation we let $(\zeta^l_\al,\eta^l_\al,g^l_\al,p^l_\al)$
denote the four components of $s^l_\al$):
\begin{enumerate}
\item  \label{clean1} For each $l<m$ the sequence $\<s^l_\al:\al\in X\>$ is tidy;
\item  \label{clean2} 
For each $\al\in X$ the sequence $\<s^l_\al:l<m\>$ is strictly $<_{I_\delta}$-increasing;
\item  \label{clean3} 
For all $l,l'<m$ and all $\al,\b\in X$,
$\zeta^l_\al<\zeta^{l'}_\al\Leftrightarrow \zeta^l_\b<\zeta^{l'}_\b$ and
$\zeta^l_\al>\zeta^{l'}_\al\Leftrightarrow \zeta^l_\b>\zeta^{l'}_\b$;
\item  \label{clean4} 
For all $l,l'<m$ and all $\al,\b\in X$,
\begin{itemize}
\item $\eta^l_\al=\eta^{l'}_\al$ if and only if $\eta^l_\b=\eta^{l'}_\b$;
\item $\eta^l_\al=\eta^{l'}_{\zeta^{l}_\al}$ if and only if $\eta^l_\b=\eta^{l'}_{\zeta^{l}_\b}$;
\item $\eta^l_\al\lessdot\eta^{l'}_{\zeta^{l}_\al}$ 
if and only if $\eta^l_\b\lessdot\eta^{l'}_{\zeta^{l}_\b}$;
\item $\eta^{l'}_{\zeta^{l}_\al}\lessdot\eta^l_\al$ if and only if 
$\eta^{l'}_{\zeta^{l}_\b}\lessdot\eta^l_\b$;
\end{itemize}
\item \label{clean5} For $\al>\b$, $\al>\zeta^l_\b$ for all $l<m$;
\item  \label{clean6} 
For all $l,l'<m$ such that $\zeta^l_\al>\zeta^{l'}_\al$ for some $\al\in X$,
$\<c(\zeta^l_\al,\zeta^{l'}_\al):\al \in X\>$ is constant;
\item  For all $l<m$ and all ordinals $\b^*$, if $\zeta^l_\b=\b^*$ for all $\b\in X$
then 
$\<c(\zeta^l_\al,\b^*):\al \in X\>$ is constant.
\end{enumerate}
}
\end{defn}

It is readily checked that if $\<b_\al:\al\in X\>$ is clean
and $b_\al=\tao(s^l_\al:l<m)$ for all $\al\in X$ then $R(s^l_\al,s^{l'}_\al)
\leftrightarrow 
R(s^l_\b,s^{l'}_\b)$ for all $l,l'<m$ and all $\al,\b\in X$.

\begin{lemma}  \label{cleanlemma}
If $S\sbq\mu$ is stationary and $\<b_\al:\al\in S\>$ is any sequence
of elements from $M_\delta$, then there is a stationary $S'\sbq S$ such that
the subsequence
$\<b_\al:\al\in S'\>$ is clean.
\end{lemma}

\bp  Since $M_\delta$ is an Ehrenfeucht-Mostowski model built from the skeleton $I_\delta$,
for each $\al\in S$ there is a term $\tao_\al$ with $m(\al)$ free variables  and
elements $s^0_\al,\dots,s^{m(\al)-1}_\al$ from $I_\delta$  such that
$b_\al=\tao_\al(s^l_\al:l<m(\al))$.  Since
$|L'|<\mu$, we can shrink $S$ to a smaller stationary set on which our choice of
$\tao$ (and hence $m$) is constant.  By applying Lemma~\ref{tidylemma} to
$\<s_\al^l:\al\in S\>$
for each $l<m$, we obtain Condition~(1).
As well, Conditions (2)--(4) and (6)--(7) are obtainable since the union of fewer than
$\mu$ non-stationary subsets of $\mu$ is non-stationary.  To obtain Condition~(5),
it suffices to note that the set
$$C=\{\al\in\mu:\al>\zeta^l_\b\ \hbox{for all}\ \b\in S\cap\al\ \hbox{and all}\ l<m\}$$ 
is club in $\mu$ (hence
$S\cap C$ is stationary).\endproof

Next we want to relate pairs of clean sequences from $M_\delta$.

\begin{defn}  \label{cleanpair}
{\em 
The (ordered) pair of sequences $\<a_\al:\al\in X\>$ and $\<b_\al:\al\in X\>$
of elements from $M_\delta$
is a {\em clean pair\/} if both sequences are clean
and the following two conditions hold (suppose that
each $a_\al=\tao_a(s^l_\al:l<m(a))$ and each 
$b_\al=\tao_b(t^{l'}_\al:l'<m(b))$):
\begin{enumerate}
\item \label{cleanpair2}
If $\zeta^{s_\al}=\al^*$ for all $\al\in X$, then $\b>\al^*$ for all $\b\in X$;
\item \label{cleanpair3}
If $\zeta^{t_\b}=\b^*$ for all $\b\in X$ then $c(\zeta^{s_\al},\b^*)=c(\zeta^{s_{\al'}},\b^*)$
for all $\al,\al'\in X$.
\end{enumerate}
}
\end{defn}

%
%
%

\begin{lemma}  \label{cleanpairlemma}
Suppose that $S\sbq\mu$ is stationary
and that
$\<a_\al:\al\in S\>$ and $\<b_\al:\al\in S\>$
are arbitrary sequences from $M_\delta$ indexed by $S$.
Then there is a stationary $S'\sbq S$ such that
the subsequences
$\<a_\al:\al\in S'\>$ and $\<b_\al:\al\in S'\>$
form a clean pair.
\end{lemma}

\bp  It follows from Lemma~\ref{cleanlemma} that we may assume that each of the
sequences is clean.  Now Condition~(1) can be obtained simply be removing
a  bounded initial segment from $S$ and Condition~(2) is obtained by noting that
there are only countably many choices for the value of $c(\zeta^{s_\al},\b^*)$
for each of the (finitely many) $\b^*$'s that are relevant.\endproof

\begin{defn}  \label{goodrelation}
{\em 
Suppose that $X\sbq\mu$.  A relation $D\sbq X^2$ {\em concentrates
on tails\/} if, for all $\al\in X$ there is $\b(\al)<\al$ such that
$$D(\al,\b)\leftrightarrow D(\al,\b')$$
for all $\b,\b'\in X$ that satisfy $\b(\al)\le\b,\b'<\al$.
}
\end{defn}

\begin{lemma}  \label{tail}  Suppose that $S\sbq\mu$ is stationary and a
relation $D\sbq S^2$ concentrates on tails.  Then there is a stationary subset
$S'\sbq S$ such that $D(\al,\b)\leftrightarrow D(\al',\b')$ for all
$\al>\b$, $\al'>\b'$ from $S'$.
\end{lemma}

\bp  Fix a function $\al\mapsto\b(\al)$ with domain $S$ that witnesses $D$ concentrating
on tails.  As this function is pressing down, it follows from Fodor's lemma that there
is a $\b^*$ and a stationary $S_1\sbq S\setminus\b^*$ such that 
$D(\al,\b)\leftrightarrow D(\al,\b')$ 
for all $\al\in S_1$ and all $\b,\b'\in S'\cap\al$.

Let $T=\{\al\in S_1: D(\al,\b)$ holds for all $\al,\b$ in $S_1$, $\al>\b\}$.
Either $T$ or $S_1\setminus T$ is stationary and hence is an appropriate choice for $S'$.\endproof

We finish this section with a type of `interpolation theorem'
for strictly increasing ordinal-valued functions.

\begin{lemma}  \label{interp}
Suppose that $S\sbq\mu$ is stationary and $\gamma$ is any ordinal.
For every strictly increasing $f:S\ra\gamma$ 
there is a club $C\sbq\mu$ and a strictly
increasing (total) function $f^*:\mu\ra\gamma$ such that
$f^*|_{S\cap C}=f|_{S\cap C}$.
\end{lemma}

\bp  First, let $B=\{\al\in S:f(\al)<f(\b)+\al$ for some $\b\in S\cap \al\}$.
If $B$ were stationary, then it would follow from Fodor's lemma that there would
be a stationary $B'\sbq B$ and a $\b^*\in S$ such that
$\al\in B'$ implies $$f(\b^*)<f(\al)<f(\b^*)+\al.$$
But then, for each $\al\in B'$
one could choose $\gamma(\al)<\al$ such that $f(\al)=f(\b^*)+\gamma(\al)$.
Another application of Fodor's lemma would show that this contradicts
the fact that $f$ is strictly increasing.  Thus, we can find a club $C_1\sbq\mu$
such that 
$f(\al)\ge f(\b)+\al$ for every pair $\al>\b$ from $S\cap C_1$.
Now define a total function $g:\mu\ra \gamma$ by:
$$g(\al)=\left\{\begin{array}{ll}
   \sup\{f(\b)+\al:\b\in S\cap\al\} & \mbox{if $S\cap C_1\cap\al\neq\emptyset$} \\
   \al                              & \mbox{if $S\cap C_1\cap\al=\emptyset$}
               \end{array}\right.$$
It is easy to verify that $C_2=\{\al\in\mu:g(\al)>g(\al')$ for all $\al'<\al\}$
is a club subset of $\mu$.  Let $S'=S\cap C_1\cap C_2$ and let $D$ be the closure of $S'$
Define a function $h:D\ra \gamma$ by:
$$h(\al)=\left\{\begin{array}{ll}
   f(\al) & \mbox{if $\al\in S'$} \\
   g(\al)                              & \mbox{if $\al\in D\setminus S'$}
               \end{array}\right.$$
It is easily checked that the function $h$ is strictly increasing on $D$.
So, let $j:\mu\ra D$ be the enumeration map (i.e., $j(\al)$ is the $\al^{\rm th}$ element of $D$)
and let $f^*:\mu\ra\gamma$ be defined by $f^*(\al)=h(j(\al))$.
The function $f^*$ is strictly increasing as both $h$ and $j$ are.
As well, the set $C_3=\{\al\in \mu:j(\al)=\al\}$ is club in $\mu$
and for $\al\in S\cap C_1\cap C_2\cap C_3$, 
$$f^*(\al)=h(j(\al))=h(\al)=f(\al)$$
so $f^*$ is as desired.\endproof

\begin{lemma}  \label{total}
Let $S\sbq \mu$ be stationary and let $g:S\ra\des{^{<\omega}}(\delta)$
be any function.  There is a stationary $S'\sbq S$ and a uniform function 
$g^*:\mu\ra
\des{^{<\omega}}(\delta)$ such that $g^*|_{S'}=g|_{S'}$.
\end{lemma}

\bp  First, by shrinking $S$ if needed, we may assume that there is a number $m$
so that $lg(g(\al))=m$ for all $\al\in S$.  Similarly, for each $i<m$ we may assume that there
is a number $n(i)$ such that $lg(g(\al))(i)=n(i)$ for all $\al\in S$.  
Let $\des_{n(i)}(\delta)$ denote the subset of $\des(\delta)$ consisting
of decreasing sequences of length $n(i)$.
Note that 
$(\des_{n(i)}(\delta),<_{lex})$ is
well ordered and hence
order-isomorphic to an ordinal.  So,  by applying Lemma~\ref{constant}
once for each $i<m$
we may assume that  each of the sequences $\<g(\al)(i):\al\in S\>$ is either $<_{lex}$-strictly increasing
or constant.  For each constant $i<m$, let $g_i$ denote its common value.
By successively applying Lemma~\ref{interp} for each non-constant $i<m$
we obtain a stationary subset $S'\sbq S$ and strictly increasing total functions $f_i:\mu\ra
\des{^{<\omega}}(\delta)$ such that $f_i(\al)=g(\al)(i)$
for all $\al\in S'$.  So define $g^*:\mu\ra\des{^{<\omega}}(\delta)$ by:
$$g^*(\al)(i)=\left\{\begin{array}{ll}
   f_i(\al) & \mbox{if $i$ is non-constant} \\
   g_i                              & \mbox{if $i$ is constant}
               \end{array}\right.$$
Clearly, $g^*$ is uniform and
$g^*|_{S'}=g|_{S'}$.\endproof
%
%

\end{document}